\documentclass{amsart}[12pt]
\usepackage{amssymb, color}

\newcommand{\C}{\mathbb{C}}
\newcommand{\Z}{\mathbb{Z}}

\newcommand{\GA}{\textup{GA}}
\newcommand{\TA}{\textup{TA}}
\newcommand{\Aff}{\textup{Aff}}
\newcommand{\lp}{\rightarrow}

\theoremstyle{definition}
\newtheorem{prob}{Problem}

\begin{document}
\begin{center}
{\bf Open problems of the conference\\
Automorphisms of Affine Spaces\\
July 6-10 2009 \\
Radboud University\\
  Nijmegen, The Netherlands}
\end{center}
(See \color{red} http://www.math.ru.nl/\%7Emaubach/AAS/index.html \color{black} for more info on this conference).\\
\ \\
{\bf Notations:} $\GA_n(\C)$ is the set of polynomial automorphisms of $\C^n$. \\
$\Aff_n(\C)$ is the set of affine automorphisms, i.e. compositions of linear maps and translations.\\
$\TA_n(\C)$ is the subset of tame polynomial automorphisms, i.e. generated by triangular and affine automorphisms.\\
\ \\

{\bf \underline{Marek Karas}}

\begin{prob} Let $F,G$ be Keller maps (i.e. polynomial endomorphisms satisfying $\det(Jac(F))=\det(Jac(G))=1$). 
Suppose
\[ F\arrowvert_{xy=0}= G\arrowvert_{xy=0} \]
Does this imply $F=G$?
\end{prob}

{\bf \underline{Vladimir Bavula}}

\begin{prob}
Is it true that $\GA_n(k)$ ($k$ a field of char. zero) is generated by tame automorphisms and finitely many one-parameter subgroups of automorphisms?
\end{prob}

{\bf \underline{Wlodzimierz Danielewski}}


\begin{prob}
Find an intrinsic algebraic invariant 
of normal finitely generated $k$-algebras, 
that would capture algebraic rigidity 
of ``tubular neighborhood of infinity'' 
and could be thought of 
as adding a ``continuous'' dimension 
to the homotopy type at infinity.
This invariant used in a straightforward way, 
although requiring perhaps complicated calculations, 
must differentiate all isomorphism classes 
of the algebras $k[x,y,z]/(x^nz-y^2-f(x)y)$.
\end{prob}

\begin{prob}
Classify normal affine equivariant embeddings 
of connected solvable linear algebraic groups 
for which all isotropy groups are semisimple.
It seems to be easy when maximal tori have dimension one.
\end{prob}

{\bf \underline{Gene Freudenburg}}

\begin{prob}
Let $R=\C[a,b]=\C^{[2]}$. Let $D:R[x,y,z]\lp R[x,y,z]$ be a locally nilpotent $R$-derivation which is triangular,
i.e. $D(x)\in R, D(y)\in R[x], D(z)\in R[x,y]$. Suppose $D$ has a slice. Does it imply that 
$\ker(D)=R^{[2]}$?\\
Note: $\ker(D)$ is an $\mathbb{A}^2$-fibration over $\mathbb{A}^2$. 
\end{prob}

{\bf \underline{Yuriy Bodnarchuk}}

\begin{prob}
Let $k$ be a field.
 Is the group $\Aff_n(k), n >2 ,$  a maximal subgroup of $\TA_n(k)$ (tame transformation group)?
\end{prob}

{\bf \underline{Eric Edo}}

\begin{prob}
Let $R:=\Z[z]/(z^3)$, $d\in R$. $P_1,Q_1\in R[y]$ such that $P_1(Q_1(y))=dy$. 
Do there exist $a,b,c\in R$ and $P,Q\in R[y]$ such that $P(Q(y))=y$ and 
$P_1(y)=aP(\frac{1}{b}y), Q_2(y)=bQ(\frac{1}{c}y)$, and $a=dc$.
\end{prob}

{\bf \underline{David Wright}}

\begin{prob}
Let $G_i$ be the subgroup of $\GA_n(k)$ that stabilises $k\oplus kx_1\oplus\ldots\oplus kx_i$.
Note that $G_n=\Aff_n(k)$.
If $n=2$ then $G_2=\Aff_2(k)$, and $G_1$ is the set of triangular automorphisms. Now $\GA_2(k)=G_1 *_{G_1\cap G_2} G_2$, the Jung-van der Kulk  theorem. \\
Question: $\GA_n(k)=<G_1,\ldots,G_n>$. And if yes, is
\[ \TA_n(k)=* G_i \]
the amalgamated product along pairwise intersections?
\end{prob}

\begin{prob}
Same question, now for the tame automorphism subgroup: define $H_i=G_i\cap \TA_n(k)$. 
Question: $\TA_n(k)=<H_1,\ldots,H_n>$. And if yes, is
\[ \TA_n(k)=* H_i \]
the amalgamated product along pairwise intersections?
\end{prob}

{\bf \underline{Wenhua Zhao}}

\begin{prob}
Let $x=(x_1, x_2, \ldots,x_n)$ and $C_n\!:=[0,1]^{\times n}$, the $n$-cube in $\mathbb R^n$. Suppose that $f(x)\in \mathbb C[x]$
and
\begin{align*}
\int_{C_n} f^m(x)\, d x=0
\end{align*}
for any $m\geq 1$.

Does this imply that $f=0$? (Note that, when $n=1$,
this is true.)
\end{prob}

A much weaker version of the problem above is
given by the next open problem. But, first let us recall the following notion.\\

{\bf Definition}\, {\it Let $R$ be any commutative ring and
$\mathcal A$ a commutative $R$-algebra.
A $R$-subspace $M$ of $\mathcal A$ is said to be
a {\it Mathieu subspace} of
$\mathcal A$  if the following property holds:
for any $a, b\in \mathcal A$ with $a^m\in M$ when
$m>>0$, we have, $a^m b \in M$ when $m>>0$.} \\

Note that, equivalently, one may replace the first ``$m>>0$" in the definition above by ``$m\ge 1$".

\begin{prob}
Let $M\!:=\{f\in \mathbb C[x] \,\,|\,\, \int_{C_n} f(x) dx =0 \}$.
Is $M$ a Mathieu subspace of the polynomial algebra
$\mathbb C[x]$?
\end{prob}

The next open problem asks if Mathieu subspaces are closed under the addition. More precisely,

\begin{prob}
Let $\mathcal A$ be a finitely generated $k$-algebra, where $k$ is any  field. Let $M_1$ and $M_2$ be any two Mathieu subspaces
of $\mathcal A$.
Is $M_1+M_2$ also a Mathieu subspace
of $\mathcal A$?
\end{prob}

For more backgrounds and motivations of the notion of Mathieu subspaces and also the open problems above,
see arXiv:0902.0212 [math.CV].\\

{\bf \underline{Leonid Makar-Limanov}}

\begin{prob}
Let $A,B$ be commutative rings. How are $ML(A), ML(B)$ and $ML(A\otimes B)$ related?
\end{prob}

\begin{prob}
Let $A,B$ be rings. Suppose that $A\otimes B = F_n$, a free associative algebra of rank $n$. Does it imply that $A,B$ are both free associative algebras too?
\end{prob}

\end{document}